\documentclass[a4paper,10pt]{amsart}

\usepackage{amstext}
\usepackage{amsmath}
\usepackage{ifthen}
\usepackage{amscd}
\usepackage{amssymb}
\usepackage[english]{babel}
\usepackage[T1]{fontenc}
\usepackage[utf8]{inputenc}
\usepackage[all]{xy}
\usepackage{graphicx}
\usepackage{enumerate}
\usepackage{xspace}
\usepackage{pdfsync}
\usepackage{epic}

\newtheorem{prop}{Proposition}[section]
\newtheorem{coro}[prop]{Corollary}
\newtheorem{lem}[prop]{Lemma}
\newtheorem{rem}[prop]{Remark}

\newtheorem{defi}[prop]{Definition}
\newtheorem{theo}[prop]{Theorem}
\newtheorem{conj}{Conjecture}

\newtheorem{ques}[prop]{Question}

\newcommand{\R}{\mathbb R}
\newcommand{\Q}{\mathbb Q}
\newcommand{\C}{\mathbb C}
\newcommand{\N}{\mathbb N}

\newcommand{\Z}{\mathbb Z}
\newcommand{\X}{\mathcal X}

\newcommand{\Od}{\mathcal{O}_d}
\newcommand{\Qd}{\Q(\sqrt{-d})}
\newcommand{\chyp}{{\mathbb H}^2_{\mathbb C}}

\newenvironment{prv}{Proof:}{$\Box$}

\begin{document}

\title[Examples of Orbifold K\"ahler Groups, II]{Orbifold K\"ahler Groups related to arithmetic complex hyperbolic lattices}

\date{April, 30 2018}
\author{Philippe {Eyssidieux}}
\thanks{This research was partially supported by the ANR project Hodgefun ANR-16-CE40-0011-01.}

\begin{abstract} We study  fundamental groups of toroidal compactifications of non compact ball quotients and show that the Shafarevich 
conjecture on holomorphic convexity for these complex projective manifolds is satisfied in dimension $2$ provided the corresponding lattice is arithmetic and small enough. 
The method is to show that  the Albanese mapping on an  \'etale covering space generates jets on the interior, if the lattice is small enough. 
We also explore some specific examples of Picard-Eisenstein type. 
\end{abstract}

\maketitle

\section{Introduction}

\subsection{}

Given $X$ a compact K\"ahler manifold, the question raised by Shafarevich whether the universal  covering space of $X$ is  holomorphically convex, also known as
the Shafarevich   Conjecture
on holomorphic convexity
$SC(X)$,  and the study of the Serre problem of characterizing the finitely presented groups arising as fundamental groups of complex algebraic manifolds 
lead to consider  several properties the fundamental group of $X$ may or may not satisfy:

\begin{enumerate}
 \item $Inf(X)$: $\pi_1(X)$ is infinite. 
 \item $Inf_{et}(X)$: Assuming $Inf(X)$, the profinite completion $\hat \pi_1(X)$ is infinite. 
 \item $RF(X)$: $\pi_1(X)$ is residually finite. 
 \item $Q(X)$: Assuming $Inf(X)$, $\pi_1(X)$ has a finite rank representation in a complex vector space whose image is infinite. 
 \item $SSC(X)$: 
For every $f:Z\to X$ where $Z$ is a compact connected complex analytic space
and $f$ is holomorphic $\#\mathrm{Im}(\pi_1(Z) \to \pi_1(X)) =+\infty \Leftrightarrow
\exists N \in \N^* \  \exists \rho: \pi_1(X) \to GL_N(\C) \  \#\rho(\mathrm{Im}(\pi_1(Z) \to \pi_1(X))=+\infty. $
 \item $L(X)$: $\pi_1(X)$ is a linear group. 
\end{enumerate}

One has the following implications, the first one being a classical result of Mal\v{c}ev, the last one was proved in increasing generality in \cite{EKPR, CCE, ajm}:
$$
L(X)\Rightarrow RF(X), SSC(X), \ \ SSC(X) \Rightarrow Q(X) \Rightarrow Inf_{et}(X),  \ \ SSC(X) \Rightarrow SC(X).
$$

The counterexamples to $RF(\_)$ -hence $L(\_)$- \cite{Tol, Trento} do not give rise to counterexamples of
$SSC(\_)$, actually all the other statements hold trivially true for the complex projective manifolds considered there. 

All these properties make sense if $X$ is replaced by a 
compact K\"ahler smooth Deligne-Mumford stack, see e.g. \cite{ajm},  a complex algebraic manifold (say quasi projective) or a smooth separated
Deligne-Mumford stack with quasi projective or quasi K\"ahler moduli space. 
Using orbifold compactifications of a given quasi-K\"ahler manifold $U$, one can produce compact K\"ahler orbifolds  with {\bf{potentially}} interesting 
fundamental groups if $\pi_1(U)$ has a sufficiently rich normal subgroup lattice. These orbifold K\"ahler groups can sometimes be proven to be fundamental groups of related 
compact K\"ahler manifolds. This happens if
 the inertia morphisms are injective after passing to the profinite completion. Checking this propery requires at least
 {\em some} understanding of the structure of the orbifold fundamental group.

A first non trivial class of $U$, complements of line arrangements in projective space, was analysed in \cite{ajm} where $SSC(\_)$ was settled affirmatively  
in the equal weight case   and it seems  much more difficult to settle $L(\_)$. 
The resulting orbifolds are abelian quotients of the Hirzebruch algebraic surfaces ramified over the arrangement \cite{Hirz, BHH}. A general theory of the fundamental group of these surfaces
in the unequal weight case seems to be out of reach except in specific cases: in
the case of  $CEVA(2)$,  which was investigated in depth for the construction of complex hyperbolic lattices,  the theory is fascinating \cite{DM,DMC,Mos}. 
 \subsection{}
The present article investigates a second non-trivial class,  where a lot of the most beautiful examples come 
precisely from the aforementioned work of Mostow and Deligne-Mostow:  finite covolume non compact quotients of a complex hyperbolic space. Then $\pi_1(U)$ is a non-uniform complex hyperbolic lattice, 
and the manifolds we investigate are the toroidal compactifications of non compact ball quotients \cite{AMRT,Mok} . These objects have attracted a lot of attention  in recent years 
from several perspectives, construction of interesting lattices \cite{DPP,Stov,DiCSto}, 
their position in the classification of algebraic varieties \cite{DiC2,BT}, Kobayashi hyperbolicity \cite{Cad1,Cad2}. 

In contrast to the case of rank $\ge 2$ where the fundamental group of a toroidal compactification of an irreducible hermitian locally symmetric space is finite, 
a toroidal compactification of a ball quotient  can have a  large fundamental group. For instance,  \cite{HS}  constructs a Riemannian metric of non positive curvature
on the toroidal compactification of a small enough lattice. These toroidal compactifications are thus $K(\pi,1)$, the universal
covering space being diffeomorphic to a real affine space and the fundamental group has exponential growth. But it does not seem possible to use the methods in \cite{HS}
to prove that the universal covering space is Stein (which what $SC$ predicts) or to construct linear representations of the fundamental group. 

\subsection{} Let us describe the content of this article. 

Given $\Gamma<PU(n,1)$ a non uniform lattice, $n\ge 2$,  decorating the construction of \cite[p. 29-30]{Hol} with its natural stack structure, or interpreting \cite[Ch. 4]{Hol3} and  \cite{Ulu}
in a more flexible language, we construct an orbifold compactifications   $[\Gamma \backslash \chyp]\subset \mathcal{X}_{\Gamma}^{tor}$
with no codimension $1$ ramification at infinity and a singular DM-stack compactification $[\Gamma \backslash \mathbb{H}^n_{\mathbb{C}}]\subset \mathcal{X}_{\Gamma}^{BBS}$. When $\Gamma$ is neat this is
the usual construction from \cite{AMRT}. 

Let $C\subset \partial  \mathbb{H}^n_{\mathbb{C}}$ be the set of the preimages of the cusps in $\Gamma \backslash  \mathbb{H}^n_{\mathbb{C}})^{BBS}$. The finite set $\Gamma \backslash C$ is the 
set of cusps of $\Gamma$. 

\begin{theo}
For each $c\in C$ of $\Gamma$ denote by $H_c$ (resp. $Z_c$) the intersection of $\Gamma$ with the unipotent radical of the parabolic subgroup attached to $c$ (resp. its center). Then:
\begin{enumerate}
\item $\pi_1( \mathcal{X}_{\Gamma}^{BBS})= \Gamma^{BBS}=\Gamma \slash < H_c, \  c \in C >$.
\item $\pi_1( \mathcal{X}_{\Gamma}^{tor})= \Gamma^{tor}=\Gamma \slash < Z_c, \  c \in C >$.
\end{enumerate}
\end{theo}

Here $<\_ >$ stands for the subgroup generated by $\_$ which is normal in the two above cases. We have not been able so far to  find any piece 
of information on these rather natural quotients of $\Gamma$ in the litterature. In spite of the fact that the credit for this result should be given to \cite{LKM} and \cite{KaSa},
we nevertheless display it in the introduction, in order  to translate the problem we study here in terms familiar to complex hyperbolic geometers. 

Then, we focus  on constructing virtually 
abelian linear representations of the fundamental group $\Gamma^{tor}$, hence on studying  the virtual first Betti number. We will give evidence for the following conjecture:
\begin{conj}
 If $\bar X_{\Gamma}$ is the toroidal compactification of a non uniform  arithmetic lattice $\Gamma<PU(n,1)$, which is torsion free and torsion free at infinity \footnote{Actually,  torsion can and will be dealt with orbifold methods.} , $SSC(\bar X_{\Gamma})$ holds and the representations can be taken 
 to be virtually abelian. Furthermore,  there is a finite \'etale Galois covering $\bar X_{\Gamma'}\to \bar X_{\Gamma}$ with Galois group $G$ such that
 the Stein factorization of the quotient Albanese morphism $\bar X_{\Gamma} \to G\backslash Alb(\bar X_{\Gamma'})$ is the Shafarevich morphism of $\bar X_{\Gamma}$.
 \end{conj}

\begin{theo}\label{thm2}
 If $\Gamma''<PU(2,1)$ is arithmetic, there is a finite index  subgroup $\Gamma'<_{fi}\Gamma''$ such that if $\Gamma<_{fi}\Gamma'$, $\bar X_{\Gamma}$ has Albanese dimension $2$ and its image contains
 no translate of an elliptic curve at its generic point. It satisfies $SC, Q$. 
 \end{theo}
 
 With the notations of Theorem  \ref{thm2},  $SSC(\bar X_{\Gamma})$ would follow from the following statement:
 $$ \leqno{(\dagger)} \quad \exists \Gamma_* \triangleleft_{fi} \Gamma'' \quad \forall c\in C \quad H_c\cap \Gamma_* \slash Z_c \cap \Gamma_* \to H_1(\Gamma_*^{tor}, \Q)
  \quad  \mathrm{is \  injective,}
 $$
 and the universal covering space of $\bar X_{\Gamma}$ would be a Stein manifold. 
 It is enough to look at all $c$ in a finite representative set for $\Gamma''\backslash C$. We write $<_{fi}$ to when we want to emphasize a subgroup inclusion has finite index.

 \begin{theo}\label{thm3}
 If $n\ge 3$ and  $\Gamma''<PU(n,1)$ is arithmetic, there is a finite index  subgroup $\Gamma'<\Gamma''$ such that if $\Gamma<_{fi}\Gamma'$, $\bar X_{\Gamma}$ has Albanese dimension $n$ and its image contains
 no translate of an abelian variety at its generic point. It satisfies $Q$. 
 \end{theo} 
 
 As a corollary of our approach, we get a new proof of the following known facts: 
 
 \begin{theo}\label{thm4}
 Under the assumptions of Theorems \ref{thm2} \ref{thm3},   $\bar X_{\Gamma}$ has ample cotangent bundle modulo its boundary, is Kobayashi hyperbolic modulo its boundary, and 
 the rational points over a  number field over which it is defined are finite modulo its boundary. 
\end{theo}

One would like effective versions of these results and quantify 
them  in terms of the ramification indices at infinity in the spirit of \cite{BT,Cad1} which  but our method is inherently non-effective.  
The examples we have studied so far, most notably the $2$-dimensional Picard-Eisenstein case where the hard work has been done by  Feustel and Holzapfel
\cite{Hol},   suggest much better statements. The article finishes with a detailled  discussion of  the Picard-Eisenstein commensurability class from the present perspective. A very strong version of $(\dagger)$ holds in this class. 

In order to streamline the discusion, we introduce in Definition \ref{refined}
  an equivalence relation, {\em  refined commensurability},  on the set of all  lattices  in a commensurability class, namely that their (orbifold) toroidal compactification are related by an 
 \'etale correspondance. There is a natural partial ordering on its 
 the quotient set, which is induced by reverse inclusion of lattices.  The questions 
 we study here depend only of the refined commensurability class and as in \cite{ajm} are  more delicate for small classes. 

\subsection{} When writing this article we were not sure whether the statement about rational points was new, it is not a corollary of \cite{Ull}.
In the final stage of the redaction,
Y. Brunebarbe informed us that,  if $n=2$,
it follows from \cite[Theorem 0.3]{Dim} a paper we were not aware of. 
The Kobayashi hyperbolicity statement is not new since it follows from \cite{Nad} and an effective, hence better, version was proved in \cite{BT} - and even
more precise results follow from \cite{Dim} if $n=2$. On the other hand \cite{Dim} does not imply $SC$.

This article  follows the same basic idea as \cite{Dim} where   N. Fakhruddin is credited for it.  
 \cite[Proposition 3.8]{Dim} gives an explicit $\Gamma$ in each 
commensurability class such that $q(X_{\Gamma})>2$. 
Here, with Lemma \ref{primitive}, we go a little further  in the study of the differential geometry of the Albanese mapping,
using as the only automorphic input the classical fact \cite{Wa} that one may achieve $q(X_{\Gamma})>0$.
We did not find a reference for the analogous property in the non arithmetic case and 
will refrain from making any conjecture in that case. 
Our method however seems to be hopelessly non-effective and  relies on the commensurator property of arithmetic lattices.

\subsection{}

These results say that   small covolume arithmetic lattices and, unsurprisingly,
non-arithmetic lattices
are  
the most interesting ones from the present perspective and we hope to come back to their study in future work.

The author would like to thank Y. Brunebarbe, B. Cadorel, B. Claudon, M. Deraux, B. Klingler, F. Sala
and G. W\"ustholz for useful conversations related to this article and the Freiburg Institute for Advanced Studies for hospitality during its preparation.

\section{Orbifold partial compactifications}

\subsection{}
As advocated by \cite{No2, art:lerman2010}, we define
an orbifold to be a smooth Deligne-Mumford  stack with trivial generic isotropy groups relative to the category of complex 
analytic spaces\footnote{Except in specific cases, one should not work  relatively to the category of complex manifolds, since the Yoneda functor  should distinguish a complex space
and the normalisation of  underlying reduced complex space.} with the classical topology: we assume that the moduli space is Hausdorff
and that the inertia groups are finite. There is an analytification $2$-functor  from DM-stacks over $\C$ to complex analytic DM-stacks and an 
underlying topological stack functor from 
complex analytic DM-stacks to topological DM-stacks \cite{No2}.

 \subsection{}
An orbifold $\X$ is said to be {\em developable} if its universal covering stack  \cite{No2} \footnote{
The covering theory of \cite{No2} has nothing 
to do with the complex structure.}
is an ordinary manifold, which is equivalent to the
injectivity of every local inertia morphism $I_x=\pi_1^{loc}(\X,x) \to \pi_1(\X,x)$, $x$ being an orbifold point of $\X$.
If $x_0$ is a base point of $\mathcal{X}$,  we have a conjugacy class
of morphisms $I_x \to \pi_1^{et}(\X, x_0)$, also called the local inertia morphisms. We will
abuse  notation and drop the base point dependency of $\pi_1$ when harmless. The orbifold $\mathcal{X}$ is said to be {\em uniformizable}
whenever the profinite completion $I_x \to \pi_1^{et}(\X)$  of every local inertia morphism is injective\footnote{The literature also uses 
{\em good} orbifold for developable and {\em very good} orbifold for uniformizable.}. When $\pi_1(\X)$ is residually finite, 
uniformizability and developability are equivalent properties. The condition  \lq residually finite\rq \  cannot be dropped, a   counterexample is given in \cite{ajm}. 
The fundamental group of a compact K\"ahler uniformizable orbifold 
occurs as the fundamental group of a compact K\"ahler manifold  \cite{vjm}.

\subsection{The fundamental group of a weighted DCN}\label{wdcn}
Let us recall the simplest  examples of orbifold compactification of quasi-K\"ahler manifolds described in \cite{ajm}.
\subsubsection{Root stacks}

Let $M$ be a (Hausdorff second countable) complex analytic space and $D$ be an effective Cartier divisor. Let $r\in \N^*$. Then, one can construct $P \to M$ the principal 
$\C^*$-bundle attached to $\mathcal{O}_M(-D)$ and the tautological section $s_D \in H^0(M,\mathcal{O}_M(D))$ can be lifted 
to a holomorphic function $f_D: P \to \C$ satisfying $f_D(p. \lambda)= \lambda f_D(p)$ for every $\lambda\in\C^*, \ p\in P$. 
Define a complex analytic space $Y:=Y_D\subset P \times \C_z$ by the equation $z^r=f_D(p)$. One can define a $\C^*$-action on $Y$ by
$(p,z). \lambda= (p\lambda^r, \lambda.z)$. The complex analytic stack $$M(\sqrt[r]{D}):=[Y_D/\C^*]$$ (see e.g. \cite{EySa, bn}) 
is then a Deligne-Mumford separated complex analytic stack with trivial generic isotropy groups whose moduli space is $M$
and an orbifold if $M$
and $D$ are smooth. The non-trivial isotropy groups live over the points of $D$ and are isomorphic to $\mu_r$
the group of $r$-roots of unity. In the smooth case, the corresponding differentiable stack can be expressed as the quotient by the natural infinitesimally free $U(1)$ action 
on the restriction $UY$ of $Y$ to the unit subbundle for (any hermitian metric) of $P$  which is a manifold indeed. It is straightforward to
see that this is an analytic version of Vistoli's root stack construction:
\begin{lem}
 If $(M, D)$ is the analytification of $(\mathcal{M}, \mathcal{D})$ a pair consisting of a $\C$-(separated) scheme and a Cartier divisor, 
 then $M(\sqrt[r]{D})$ is the analytification of $\mathcal{M}_{O (\mathcal{D}), s_{\mathcal{D}}, r}$ in the notation of \cite[Section 2]{cadman2007}. 
\end{lem}

A Cartier divisor $D$ on a scheme $M$ should be thought of as a pair of an invertible sheaf $\mathcal{L}$ and a section $s_{D} : \mathcal{O} \to \mathcal{L}$
such that $[s_D]=D$. In other words a map of algebraic stacks $\mu: M \to [\C \slash \C^*]$. We have the natural $r$-th power 
map $.^r:[\C \slash \C^*] \to [\C \slash \C^*]$ and an equivalence $M(\sqrt[r]{D}) \to M\times_{[\C \slash \C^*] \mu, .^r}  [\C \slash \C^*]$. 
Using this we may even promote $M$ to be a stack and allow for $s_{D}=0$. In the latter case we get a $\Z\slash r\Z$-gerbe on $M$
whose class is the reduction mod $r$ of $c_1(\mathcal{L})$.
 The main property  is treated in the scheme-theoretic setting by \cite{cadman2007}:  
 \begin{lem}\label{morphstack} If $S$ is a complex analytic stack 
 $$Hom(S,M(\sqrt[r]{D})) = \{f: S\to M \ \mathrm{\C-analytic}, \exists D_S \  \mathrm{Cartier \ on \ } S \   \mathrm{s. t.} \ D_S=r.f^*D\}. $$
 \end{lem}

\subsubsection{}

Let $\bar X$ be a compact K\"ahler manifold and $x_0\in X$ a base point, $n=\dim_{\C}(X)$, 
and let $D:=D_1+ \ldots + D_l$ be a simple normal crossing divisor whose smooth irreducible components are denoted by $D_i$. We assume for simplicity $x_0\not \in D$. 
For each choice of weights $d:=(d_1, \ldots, d_l)$, $d_i\in \N^*$,  one may construct 
as \cite[Definition 2.2.4]{cadman2007} does in the setting of scheme theory 
the compact K\"ahler orbifold (Compact K\"ahler DM stack with trivial generic isotropy) 
$$\X(\bar X, D, d):= \bar X(\sqrt[d_1]{D_1} )\times_{\bar X} \ldots \times_{\bar X} \bar X(\sqrt[d_l]{D_l}).$$
In other words, $\X(\bar X, D, d)=[Y_{D_1}\times_{\bar X} \dots \times_{\bar X} Y_{D_l} / \C^{*l}]$. 
Denote by $X$ the quasi-K\"ahler manifold $X:=\bar X \setminus D$. View $\X(\bar X, D, d)$ as an orbifold compactification of $X$ and denote by $j_d: X \hookrightarrow \X(\bar X, D, d)$ the natural open immersion. 

By Zariski-Van Kampen, 
the fundamental group $\pi_1(\bar X,x_0)$ is the quotient of $\pi_1(X,x_0)$ by the normal subgroup generated by the $\gamma_i$, where $\gamma_i$ is a meridian loop for $D_i$. 
Zariski-Van Kampen generalizes to orbifolds, see e.g. \cite{No1,Zoo}, and $\pi_1(\X(\bar X, D, d),x_0)$ is the quotient of $\pi_1(X,x_0)$ by the normal 
subgroup generated by the $\gamma_i^{d_i}$. 

\begin{rem}
If $(X,\Delta)$ is an orbifolde in the sense of Campana, the Campana orbifold fundamental group is the fundamental group
of the root stack on $X\setminus D^1$ where $D^1$ is the log-singular set of $(X,\Delta)$. It may happens that the fundamental group 
of the trace of the root stack on $X\setminus D^1$ of a small ball centered on $D^1$ is finite and in this case we get a DM-stack compactification which may be singular. 
In dimension 2, the list of local configurations giving rise to an orbifold compactification with a  smooth moduli space is given in \cite{Ulu}. 
\end{rem}

\begin{rem}
 The orbifolds constructed here are specified up to equivalence by their moduli space and ramification indices in codimension 1 by \cite{GerSa}. 
 Hence, in dimension 2, these orbifolds carry the same information  as Holzapfel's orbital surfaces or Uludag's orbifaces. One could also use orbifolds in the sense 
 of Thurston in higher dimension. However, it is convenient to have
 at our disposal 
 maps of stacks (defined as functors of fibered categories), moduli spaces, substacks, (2-)fiber products hence fibers, 
 basic homotopy theory \cite{No2} and many differential  geometric constructions \cite{EySa}. 
\end{rem}

\section{Toroidal compactifications of complex hyperbolic orbisurfaces} 

\subsection{Orbifold Toroidal compactifications} 
\subsubsection{} The following constructions are well known to the experts \cite{BT} but we recall them in order to fix the notations.
Let $\chyp:=PU(1,2)\slash P(U(1)\times U(2))$ be the complex hyperbolic plane. 
 Let $\Gamma < PU(1,2)$ be non uniform lattice. Then $X_{\Gamma}:=\Gamma \backslash \chyp$ is an algebraic 
 quasi-projective variety by Baily-Borel's theorem in the arithmetic case, 
and \cite{Mok} for the non arithmetic case. It has a
complex projective compactification $X_{\Gamma} \subset X_{\Gamma}^{BBS}$ obtained by adding a finite number of points we shall refer to as cusps. 
 
Each cusp $c$ has some neighborhood $V_c$ such that the preimage  of $U_c:=V_c\setminus\{ c \}$ is a disjoint union of horoballs $W_{\tilde c}$ of $\chyp$ 
labelled by $ \partial \chyp \ni \tilde c  := \partial \chyp \cap \overline{ W_{\tilde c}}$ where the closure  and boundary are taken in the euclidean topology
of the projective plane $\mathbb{P}$  dual to $\chyp$. The set $C$ of all such $\tilde c$ is acted upon with finitely many orbits by $\Gamma$. 
For such a $\tilde c$, let $\Gamma_{\tilde c}$
be the stabilizer of $\tilde c$ in $\Gamma$. 

Since the stabilizer $S_{\tilde c}$ of $\tilde c$ in $PU(1,2)$ has a 3 dimensional
Heisenberg group of unipotent $3\times 3$ upper triangular  matrices with real coefficients
as its uniportent radical $H_{\tilde c}$ and $\C^*$ has its Levi component,  see \cite{Par}, we have 
an exact sequence:
$$
1 \to H_{\tilde c} \cap \Gamma_{\tilde c} \to \Gamma_{\tilde c} \to \mu_{k_c} \to 1
$$
where $\mu_k < \C^*$ is group of $k_c$-roots of unity. Whenever this causes no confusion,  we use a shorthand notation $k=k_c$.
Then, $ H_{\tilde c} \cap \Gamma_{\tilde c} <  H_{\tilde c} $ is a lattice in $H_{\tilde c}$
One has $\gamma W_{\tilde c} \cap W_{\tilde c} \not = \emptyset \Leftrightarrow \gamma \in \Gamma_{\tilde c}$. 
The center $Z_{\tilde c}$ of $H_{\tilde c} \cap \Gamma_{\tilde c}$ is a cyclic infinite subgroup and we have $Z_{\tilde c}=Z(H_{\tilde c})\cap \Gamma_{\tilde c}$. 

The group  $\Lambda_{\tilde c}:= Z_{\tilde c} \backslash H_{\tilde c} \cap \Gamma_{\tilde c}  $ is a lattice in the real 2-dimensional 
additive group $A_{\tilde c}^\tau:=Z(H_{\tilde c}) \backslash H_{\tilde c}$ 
which has a natural  structure of an affine complex line $A_{\tilde c}$. Hence $A^\tau_{\tilde c}$ is naturally a one dimensionnal complex additive group
acting as the translation group on $A_{\tilde c}$. This complex structure comes 
from the fact that $Z(H_{\tilde c})$ stabilizes a unique complex geodesic having  
$\tilde c$ as a boundary point. All such complex geodesics are the trace of a complex projective line  in $\mathbb{P}$ through $\tilde c$
and form a single orbit of $H_{\tilde c}$. In particular we have a bijection
$ A_{\tilde c}^\tau \to \mathbb{P}(T_p\mathbb{P}) \setminus \{l_{\tilde c}\}$ 
where $l_{\tilde c}$ is the complex line tangent to $\partial \chyp$ at $\tilde c$.
The linear projection from $\tilde c$ give an equivariant holomorphic map $W_{\tilde c} \to  A_{\tilde c}$  where $S_{\tilde c}$ acts through its quotient group
$Z_{\tilde c} \backslash S_{\tilde c} \simeq A^{\tau}_{\tilde c} \rtimes \C^*$. The latter group acts as the complex affine group of $A_{\tilde c}$. 

The natural map $\psi:\tilde V^1_{\tilde c}:=Z_{\tilde c}\backslash W_{\tilde c} \to A_{\tilde c}$ is a holomorphic fiber  bundle whose fiber at $\lambda \in A_{\tilde c}$ is
$Z_{\tilde c} \backslash\lambda \cap W_{\tilde c}$ a pointed disk. There is an effective action of $Z_{\tilde c} \backslash S_{\tilde c}$ on $\tilde V^1_{\tilde c}$
hence an action of $A_{\tilde c}^{\tau}$ such that $\psi$ is equivariant. 
Now, consider the genus one Riemann Surface  $E_{\tilde c}:=\Lambda_{\tilde c} \backslash A_{\tilde c}$.  The commutator gives a natural symplectic form 
on $\Lambda_{\tilde c}$ with values in $Z_{\tilde c}$ hence, with an adequate choice of a generator of $Z_{\tilde c}$,  a polarization
$\Theta_{\tilde c}$ of the weight $-1$ Hodge structure on $\Lambda_c$ induced by 
the complex structure on $\Lambda_{\tilde c} \otimes_{\Z}\R$.  
There is also a map $\psi': V^1_{\tilde c}:=\Lambda_{\tilde c} \backslash \tilde  V^1_{\tilde c} \to E_{\tilde c}$ which is a holomorphic fiber bundle
in pointed disks. A coordinate calculation enables to see that it is biholomorphic to the complement of the zero section 
of a unit disk bundle of a hermitian line bundle $(L_{\tilde c}, h_{\tilde c})$ on  $E_{\tilde c}$ with constant curvature 
 $c_1(L_{\tilde c})=-\Theta_{\tilde c}$.  The degree of $\Theta_{\tilde c}$ is the index in $Z_{\tilde c}$ of the image of the symplectic form on $\Lambda_{\tilde c}$. 
 
Hence there is a partial compactification $V^1_{\tilde c} \simeq U^{1, tor}_{\tilde c} \setminus D^{1, tor}_{\tilde c}$ where 
$U^{1, tor}_{\tilde c}$ is the full unit bundle of $(L_{\tilde c}, h_{\tilde c})$ and $D^{1, tor}_{\tilde c}$ is the zero section, a smooth divisor 
isomorphic to $E_{\tilde c}$. 

The quotient action of the quotient a group $\mu_k$ on the hermitian line bundle $(L_{\tilde c},h_{\tilde c})$ gives an action on $U^{1, tor}_{\tilde c}$ in such a way that the natural 
retraction $\pi^1_{\tilde c}: (U^{1, tor}_{\tilde c}, D^{1,tor}_{\tilde c}) \to  E_{\tilde c}$ is equivariant. In particular the group $\mu_k$
acts on $E_{\tilde c}$ in such a way that the action on $H^0(E_{\tilde c}, \Omega^1)$ is given by complex multiplication. Since it is an action 
by automorphisms and since the Lefschetz number of an automorphism is an integer it follows that $k \in \{ 1, 2, 3,4, 6 \}$\footnote{A torsion free lattice will be called neat if $k=1$ for every cusp.}. 

Dividing out $\pi^1_{\tilde c}$ by $\mu_k$ we get a retraction $\pi^1_{\tilde c}: U^{tor}_{\tilde c} \to D^{tor}_{\tilde c}$ 
in such a way that $U_c \simeq U^{tor}_{\tilde c} \setminus D^{tor}_{\tilde c}$ and a map $U_{\tilde c}^{tor} \to U_{c}$ contracting
$D^{tor}_{\tilde c}$ to $\tilde c$.

Since the construction is indepedant of $\tilde c$ we may drop this dependency in our notation and redefine $(U^{1,tor}_c, D^{1, tor}_c,E_c):= (U^{1,tor}_{\tilde c}, 
D^{1, tor}_{\tilde c},E_{\tilde c}) $ for some $\tilde c$. 

Gluing the $U^{tor}_{c}$ with $X_{\Gamma}$ along $U_c$ gives a normal surface $X^{tor}_{\Gamma}$ with a family $(D_c)_{c\in C}$ of disjoint curves with an isomorphism
$X_{\Gamma} \to X^{tor}_{\Gamma} \setminus \bigcup_{c\in C} D_c$ and a map $X^{tor}_{\Gamma} \to X^{BBS}_{\Gamma}$ contracting $D_c$ to $c$. 
Then, $X^{tor}_{\Gamma}$ is a normal surface with quotient singularities and is projective algebraic too \cite{AMRT,Mok}. 

It is actually better to glue the orbifold $ \mathcal{X}_{\Gamma}:=[\Gamma \backslash \chyp]$ (stack theoretic quotient)
with $[\mu_k \backslash U^{1, tor}_{\tilde c}]$ along $[\Gamma_{\tilde c} \backslash W_{\tilde c}]$ to get
a compact orbifold $\mathcal{X}_{\Gamma}^{tor}$ containing $ \mathcal{X}_{\Gamma}$ as the complement of a smooth divisor $\mathcal{D}_{\Gamma}$ consisting of a finite number 
of disjoint smooth substacks $(\mathcal{D}_c)_{c\in C}$ whose generic points have no inertia. One also has a stack theoretic compactification $\mathcal{X}_{\Gamma}^{BBS}$ obtained by adding 
the disjoint union of the $(B\mu_{k_c})_{c\in C}$ and a contraction map $\mathcal{X}_{\Gamma}^{tor} \to \mathcal{X}_{\Gamma}^{BBS}$. 

When the lattice is neat the stack we constructed
is equivalent to the usual smooth toroidal compactification of \cite{AMRT,Mok}. On the other other hand, $\mathcal{X}_{\Gamma}^{tor}$ is NOT the
quotient stack of the toroidal compactification attached to a neat normal sublattice: the generic point of the boundary has trivial isotropy.

\subsubsection{}

Some comments have to be made regarding the gluing construction we are performing. First of all, gluing DM topological stacks along open substacks 
is always possible thanks to \cite[Cor. 16.11, p. 57]{No2} - we are using the class of local homeomorphisms as  {\bf{LF}}.  In order to have a better picture of the toroidal compactifications, they are smooth complex DM-stacks, we can present the open substacks as the quotient 
of their frame bundles by the general linear group. The frame bundle is indeed representable, e.g. an ordinary complex manifold, and carries an {\em infinitesimally} free
proper action of the general linear group. If this action is free the stack is equivalent to an ordinary manifold. An equivalence of smooth DM stacks then 
gives rise to an isomorphism 
of the frame bundles intertwining the infinitesimally free actions and these glue perfectly well along invariant open subsets. 

Also, the construction can be performed with a non-effective finite kernel action whose image is a lattice the price being that one has to 
consider general smooth Deligne-Mumford stacks. This seems to be inevitable if one wants to work with lattices in $U(2,1)$ as in \cite{Hol}.

\subsection{Fundamental Groups of orbifold toroidal compactifications}

The moduli space of $\mathcal{X}_{\Gamma}^{*}$  is $X_{\Gamma}^{*}$ for $*=\emptyset , BBS, tor$. 
We will denote by $m:\mathcal{X}_{\Gamma}^{*}\to X_{\Gamma}^{*}$ the moduli map.   Using Van Kampen \cite{Zoo}, we get:

\begin{lem}\label{vk}

Let $x_0$ be a base point of $\mathcal{X}_{\Gamma}$. Then $\pi_1(\mathcal{X}_{\Gamma},x_0)=\Gamma$, 
\begin{itemize}
 \item 
$\pi_1(\mathcal{X}^{tor}_{\Gamma},x_0)$ is the quotient $\Gamma^{tor}$ of $\Gamma$ by the
subgroup normally generated by the $Z_{c}$
\item 
$\pi_1(\mathcal{X}^{BBS}_{\Gamma},x_0)$ is the quotient $\Gamma^{BBS}$ of $\Gamma$ by the
subgroup normally generated by the $H_{ c}\cap \Gamma_{ c}$.
\end{itemize}

\end{lem}

\begin{coro}[\cite{KMR}]The natural map
 $H^1(\Gamma^{tor}, \Q) \to H^1(\Gamma, \Q)$ is an isomorphism, consequently the Deligne MHS on $H^1(\Gamma \backslash \chyp, \Z)$ is pure of weight one. 
\end{coro}
\begin{prv}
Dually, it is enough to show that $H_1(\Gamma^{tor}, \Q) \leftarrow H_1(\Gamma, \Q)$ is an isomorphism. 
This amounts to proving that the image of $H_1(Z_{ c}) \to H_1(\Gamma)$ is torsion. 
The group $Z_{ c}$ contains the commutator subgroup of $H_{ c} \cap \Gamma_{ c}$ as a finite index subgroup. In particular it maps to the torsion
subgroup of $H_1(\Gamma)=\Gamma \slash [\Gamma, \Gamma]$. 
\end{prv}

The moduli map $m$ presents the fundamental group of $X_{\Gamma}^{*}$  is the quotient of the fundamental group of $\mathcal{X}^{*}_{\Gamma}$ by the normal subgroup 
generated by the images of the inertia morphisms, see \cite{No1}, and gives an isomorphism on cohomology with rational coefficients. 

\begin{rem} Lemma \ref{vk} is not new, the first point is the easiest special case of \cite{KaSa}, the second point results from \cite{LKM}. Note that in the  rank $\ge 2$ case, Margulis' normal subgroup theorem 
implies that the fundamental group of a toroidal compactification of a $\R$-rank $\ge 2$ irreducible locally hermitian symmetric space 
is finite. 
\end{rem}

\begin{ques}
 Can $\Gamma^{tor}$ be finite? 
\end{ques}

\subsection{Ramification of the natural  Orbifold Toroidal compactifications maps}

\begin{lem}\label{ramif}
Let $\Gamma'< \Gamma$ be a finite index subgroup. Then,  the finite covering map $X_{\Gamma'} \to X_{\Gamma}$ lifts to an orbifold map
$\mathcal{X}^{tor}_{\Gamma'} \to \mathcal{X}^{tor}_{\Gamma'}$ which restricts over $X_{\Gamma}$ to an \'etale map $\mathcal{X}_{\Gamma'} \to \mathcal{X}_{\Gamma}$.

\begin{itemize}
 \item 
Let $c'$, $c$ be cusps such that the mapping $\eta: X^{BBS}_{\Gamma'} \to X^{BBS}_{\Gamma}$ 
maps $c'$ to
$c$. The ramification index of $\mathcal{D}_{c'}$ over $\mathcal{D}_c$ is $d_{c',c}:=[Z_{c'}:Z_{c}]$. 
\item If $\Gamma'$ is normal in $\Gamma$ and $G=\Gamma' \backslash \Gamma$ is the quotient subgroup,  $G$ acts \footnote{in the sense of \cite{Rom}.
Actually $G$ acts on the frame bundle of $\mathcal{X}^{tor}_{\Gamma'}$, an ordinary complex manifold, and 
the action commutes with the right action of the general linear group, from this it is easy to find an \'etale chart
with a strict action of $G$ on the corresponding \'etale groupoid.} on $\mathcal{X}^{tor}_{\Gamma'}$,
$d_{c',c}:=d_{c}$ depends only on $c$
and $[G\backslash \mathcal{X}^{tor}_{\Gamma'}]= \mathcal{X}^{tor}_{\Gamma}(\mathcal{D}_{\Gamma}, (d_c)_{c\in C})$. 
\end{itemize}

\end{lem}

Due to ramification,  the $\Gamma^{tor}, \  \Gamma \in \mathcal{C}_d$,   split into infinitely many commensurability classes.

\subsection{Hermitian forms over $\Q(\sqrt{-d})$}

Let $d\in \N^*$ be a squarefree positive integer. 

The imaginary quadratic field $\Q(\sqrt{-d})$  has the complex conjugation as its Galois isomorphism. A non degenerate hermitian form $H$ over $\Qd$
defines a $\Q$-algebraic group $U(H)$ which is a form of $GL(\dim H)$. The $\Qd$ vector space 
$V_H$ underlying $H$ will be denoted by $V_H$. The signature of $H$ is the signature of the corresponding hermitian form also denoted by $H$   on $V:=V_H\otimes_ {Q} R:=H\otimes_{\Qd} \C$ where $\Qd\to \C$ is a given embedding.  

Let $\Od \subset \Q(\sqrt{-d})$ be the subring of its quadratic integers. Let $L$ be a free $\Od$ module which is a lattice
in the $\Qd$ vector space 
$V_H$ underlying $H$. Then $\Gamma_{H, L} = PU(H)\cap PAut(L)$ is an arithmetic subgroup and the $\Gamma_{H,L}$ all belong to a commensurability class $\mathcal{C}_d$ of
lattices in $PU(H)$. Denote by $\Gamma_{\Q}$ be the group of $\Q$-points of $PU(H)$ which is dense in $PU(H_{\R})\simeq PU(2,1)$ with respect to the classical topology. 

The study of the groups $\pi_1(\mathcal{X}^{tor}_{\Gamma})$ and $\pi_1(\mathcal{X}^{tor}_{BBS})$ for $\Gamma\in \mathcal{C}_d$
does not seem to have been carried out systematically in the litterature even in that simple case.

It is known that all non uniform commensurability classes of arithmetic lattices in $PU(2,1)$ are of the form $\mathcal{C}_d$. 

\section{Proof of the main theorems}

\subsection{Commensurability classes with non vanishing virtual $b_1$}

By commensurability, we mean commensurability in the wide sense: 

\begin{defi} Two lattices $\Gamma_1, \Gamma_2<PU(2,1)$ are commensurable if there exists a finite index torsion free lattice $\Gamma'_i < \Gamma_i$ and a holomorphic
isometry $\Gamma'_1\backslash \chyp \to \Gamma'_2\backslash \chyp $. 
 \end{defi}

\begin{defi}
 Say a commensurability class $\mathcal{C}$ of non uniform lattices in $PU(2,1)$ has  non vanisihing virtual $b_1$ if some member $\Gamma_0$ satisfies $b_1(\Gamma_0, \Q)>0$.
\end{defi}

Each commensurability class of non uniform arithmetic lattices in $PU(2,1)$ has  non vanishing virtual $b_1$ \cite{Wa}, generalizing \cite{Ka}. 
Since these lattices are finitely generated and passing to a sublattice increases $b_1(\_, \Q)$:

\begin{lem}
 If $\mathcal{C}$  has  non vanisihing virtual $b_1$ every member $\Gamma$ has a finite index normal subgroup $\Gamma_1$ such that $b_1(\Gamma_1, \Q)>0$.
\end{lem}

\subsection{Freeness of virtual Albanese Mappings}

We will now fix an {\bf{arithmetic}}  lattice which is torsion free and has unipotent mondromies so that the toroidal compatification $\bar X$ is a complex projective manifold.
Using $b_1(\Gamma, \Q)= b_1(\Gamma^{tor}, \Q)$ we conclude that there is a non zero closed holomorphic one-form $\alpha$ on $\bar X$. We restrict $\alpha$ to $X$ 
lift it to $\chyp$ to construct a closed one form $\tilde \omega \in \Omega^1(\chyp)$. 

\begin{lem}
 Every element of the vector space $V$ spanned by the $\Gamma_{\Q}.\tilde \omega$ is the lift of a holomorphic closed one form on $X_{\Gamma'}^{tor}$ for some  normal finite index sugroup
 $\Gamma'\le \Gamma$. 
\end{lem}
\begin{prv}
 Indeed $\Gamma_{\Q}$ is the commensurator of $\Gamma$. Here we use arithmeticity in a crucial way.
\end{prv}

\begin{lem}
 The closure $\bar V$ in the Fr\'echet space $\Omega^1(\chyp)$ is a non zero vector space of closed holomorphic $1$-forms which is preserved by the action of $PU(H_{\R})$. 
\end{lem}
\begin{prv}
 Indeed $\Gamma_{\Q}$ is dense in  $PU(H_{\R})$.
\end{prv}

\begin{coro}\label{coframe}
 For every  point $o$ of $\chyp$ there are two elements of $V$ which  gives a coframe of the tangent space at $o$. 
\end{coro}
\begin{prv}
 The restriction $r:\bar V \to \Omega^1_{\chyp, o}$ is equivariant under the stabilizer $K$ of $o$. Hence the image of this linear map is $K$-equivariant and no zero, the 
 surjectivity of $r$ follows from the irreducibility of the isotropy action on the
 cotangent space. 
\end{prv}

\begin{coro} \label{unram}
 There is a finite index normal subgroup $\Gamma'$ such that every $\Gamma''<\Gamma'$,  has an Albanese map which is unramified on $X_{\Gamma''}$. 
 In particular  $\bar{X}^{tor}_{\Gamma''}$ satisfies the Shafarevich conjecture.
\end{coro}

\begin{prv}  The 
immersivity outside the boundary is an immediate consequence of Corollary \ref{coframe} using noetherian induction. 
I learned from \cite{Stov}  that the fact that the Albanese map does not factor through a curve follows from \cite{Clo}.
The application to Shafarevich conjecture is then a consequence of \cite{Nap} since an {\bf{irreducible}} 
connected component of a fibre of the Albanese mapping cannot give rise to a Nori chain. Note that $SSC$ need not be satisfied. 
\end{prv}

\begin{lem} \label{primitive} The space $W=\int \bar V \subset \mathcal{O}
 (\chyp)$ of primitives of the elements of  $\bar V$ is infinite dimensional. Actually the map $W \to \mathcal{O}_{\chyp,o}/\mathfrak{m}^N$ is surjective for all $N>0$. 
\end{lem}

\begin{prv} $W$  is invariant under the whole group $PU(H_{\R})$. The evaluation map $ev:W \to \mathcal{O}_{\chyp,o}$ is $K$-equivariant.
So is the completed evaluation map: $W \to \widehat{\mathcal{O}}_{\chyp,o}$. But this is $K$ equivariantly isomorphic to
$\C[[\mathfrak{m}\slash\mathfrak{m}^2]]=\sum_{n\in \N} Sym^n \mathfrak{m}\slash\mathfrak{m}^2$. To construct this isomorphism we have chosen the polynomials as a $K$-invariant
subspace of holomorphic functions and the linear functions as generators of the maximal ideal  $\mathfrak{m}$. The $K$ equivariance implies
$ev(W)\supset ev(W)\cap Sym^n \mathfrak{m}\slash\mathfrak{m}^2$ for all $n$. If $W$ were finite there would be a finite number of nontrivial
$ev(W)\cap Sym^n \mathfrak{m}\slash\mathfrak{m}^2$ which would give a direct sum decomposition of $ev(W)$. Hence $W$ would consist 
of polynomials. Restricting to a generic line through the origin (viewed as a complex geodesic curve isomorphic to a unit disk) there would be 
a non zero space of polynomials in one variable,  generated by monomials,
fixed by the homographies which are the automorphisms of the unit disk, which is obviously absurd.

We have
already proved surjectivity of $\bar V \to \mathcal{O}_{\chyp,o}/\mathfrak{m}^2$. We now use induction on $N$ and assume the theorem is proved for some $N$.
Consider the smallest integer $M>N$ such that $ev(W)\cap Sym^M \mathfrak{m}\slash\mathfrak{m}^2\not =0$. $M$ exists since 
$W$ is infinite dimensional. Using the irreducibilty of the isotropy representation  on $Sym^M \mathfrak{m}\slash\mathfrak{m}^2$ we see that 
$ev(W)\cap Sym^M \mathfrak{m}\slash\mathfrak{m}^2=Sym^M \mathfrak{m}\slash\mathfrak{m}^2$. Hence there is a monomial say $z_1^M$ in 
$W$. Now one of the infinitesimal generators of the Lie algebra of $PU(H_{\R})$ takes the form $\xi=\frac{\partial}{\partial z_1} + \xi'$ where 
$\xi'\in \mathfrak{m} \mathcal{T}_{\chyp}$. If $M>N+1$,  $\xi . z_1^M \in W$ would contradict the minimality of $M$. Hence $M=N+1$ which is the desired conclusion. 
\end{prv}

\begin{coro} Fix $M\in \N$. 
 There is a finite index normal subgroup $\Gamma'$ such that forall $\Gamma''<\Gamma'$,  $\bar{X}^{tor}_{\Gamma''}$ has an Albanese map which generates 
 $M$ jets at all points of $X_{\Gamma''}$. 
\end{coro}

 In particular, if $M=2$ the genus zero curves in the Albanese image lie on the image of the boundary and we recover special cases of classical results:
 \begin{coro} [\cite{Nad, BT}]
$\bar{X}^{tor}_{\Gamma''}$
 satisfies the Green-Griffiths-Lang conjecture
  \end{coro}  
  \begin{prv}
 Immediate corollary of \cite[Theorem D]{Och}. 
  \end{prv}

  \begin{coro}
  $X_{\Gamma''}$ has finitely many rational points on any number field of definition.
 \end{coro}
 \begin{prv}
  Immediate corollary of Faltings' solution of Lang's conjecture \cite{Fal}, see also \cite[Theorem F.1.1, p. 436]{HinSil}. 
 \end{prv}

\subsection{Refined commensurability classes} 

The following definition is a slight generalization of the notion of completely \'etale map between negatively elliptic bounded surfaces \cite[p.256]{Hol2}:

\begin{defi}\label{refined} Two (commensurable) lattices $\Gamma_1, \Gamma_2<PU(2,1)$ are refined commensurable if there exists a finite index lattice $\Gamma'_i < \Gamma_i$ and a holomorphic
isometry $\Gamma'_1\backslash \chyp \to \Gamma'_2\backslash \chyp $ such that $\mathcal{X}_{\Gamma'_i}^{tor} \to \mathcal{X}_{\Gamma_i}^{tor}$ is \'etale. 
 \end{defi}

 If $ \mathcal{X}_{\Gamma_i}^{tor}$ is uniformizable we can assume that $\Gamma'_i$ is torsion free and neat. The questions (1)-(6) (and the orbifold Kodaira dimension) in the introduction depend only
 on the refined commensurability class. 
 
 \begin{lem}
 Two commensurable arithmetic lattices $\Gamma_1,\Gamma_2$ having a common finite index lattice are refined commensurable iff for every $\tilde c \in C \subset\partial \chyp$
 there is a common finite index subgroup $\Gamma_3$ such that  $Z(H_{\tilde c})\cap \Gamma_1=Z(H_{\tilde c})\cap \Gamma_2=Z(H_{\tilde c})\cap \Gamma_3$. 
 \end{lem}

There is an order on the set $\mathcal{C}_r$ of refined commensurability classes of a given commensurability class. We say $\mathcal{C}_1 \prec \mathcal{C}_2$ is there is some member 
of $ \mathcal{C}_2$ is a finite index subgroup of a member of  $ \mathcal{C}_1$.  It is clear that $(\mathcal{C}_r,\prec)$ is a filtering order and that the questions investigated 
in the introduction are more difficult for small elements (with the exception of $L(\_)$).
\begin{ques}
 Is $(\mathcal{C}_r,\prec)$ Artinian? Is every initial segment finite? 
\end{ques}

\begin{defi}
 A refined commensurability class is said to be small if the common universal covering stack of the corresponding $\mathcal{X}^{tor}_{\Gamma}$ is not a Stein manifold. 
 Its $\Gamma$-dimension is the 
 dimension of the Campana quotient of the moduli space of the universal covering stack by its compact complex subvarieties. 
\end{defi}

They are the most interesting classes from the present perspective.

\subsection{Higher dimensions}

The only missing ingredient being \cite{Nap}, we get if $n\ge 3$: 
\begin{theo}
 If $\Gamma<PU(n,1)$ is arithmetic and  small enough in its commensurability class, $\bar X_{\Gamma}$ has Albanese dimension $n$ and its image contains
 no translate of an abelian subvariety at its generic point. It satisfies $ Q$, 
 has big  cotangent bundle, is Kobayashi hyperbolic modulo its boundary, and 
 the rational points over a  number field over which it is defined are finite modulo its boundary. 
\end{theo}

\section{The Picard-Eisenstein commensurability class}

\subsection{The Picard-Eisenstein Lattice}
Let us look at the Picard-Eisenstein lattices using its beautiful two-generator  presentation \cite{FP}. This group, denoted by $PU(H_0, \Z [\omega])$ 
is in $\mathcal{C}_3$ and is actually of the form  $\Gamma_{H, L}$ for the hermitian form whose matrix is $H_0=
\left(
\begin{array}{ccc}
 0&0&1 \\
 0&1&0 \\
 1&0&0
\end{array}
\right)
$ and $L=\mathcal {O}_3^{\oplus 3}=\Z[\frac{-1+i\sqrt{3}}{2}]^{\oplus 3}$ is the standard lattice. 

\begin{prop}\label{surj}
$PU(H_0, \Z [\omega])^{tor}$ has a surjective morphism to a $(2,3,6)$ orbifold group.
\end{prop}

\begin{prv}
 The presentation in \cite{FP} is
 $$ <P,Q,R| R^2=[PQ^{-1} , R]= (RP)^3=P^3Q^{-2}=(QP^{-1})^6>
 .$$ We are imposing $P^3=Q^2=1$ in this presentation by \cite[p. 258]{FP}. Further imposing $R=1$ gives $<P,Q| P^3=Q^2=(PQ^{-1})^6=1>$ as a quotient group. 
\end{prv}

It may be reassuring to get a confirmation of the:

\begin{coro}
 The Picard-Eisenstein commensurability class $\mathcal{C}_3$ has  non vanishing virtual $b_1$.
\end{coro}

\subsection{Hirzebruch's example of a configuration of elliptic curves in an abelian surface whose  complement is complex hyperbolic} \label{sec2}

It is a well known theorem of Holzapfel \cite{Hol2} that the following  example is in  $\mathcal{C}_3$  .

\subsubsection{}
Let $E\simeq \C / \Z[j]$ be the elliptic curve isomorphic to Fermat's cubic curve and let $\bar X$ be the blow up of $(0_{E \times E})$. We denote by $D_i$ the strict transforms in $\bar X$
of the  ellitic curves $E_i$ where:
$$E_1=E\times \{ 0_E\}, E_2= \{ 0_E\}\times E, E_3=\Delta=Gr(id_E), E_4=Gr(-j)$$
and we  obtain a SNC divisor $(\bar X, D)$ where $D=D_1+D_2+D_3+D_4$. The $D_i$ are pairwise disjoint. 
\begin{theo}(Hirzebruch, \cite{BHH}) $X\simeq \Gamma_{Hirz} \backslash \chyp$ where $\Gamma_{Hirz} \subset PU(2,1)$ is a non uniform lattice. 
$\bar X$ is its toroidal compactification. 
 \end{theo}
 
 The notation $(\bar X, D)$ will be used throughout this section to denote this particular confinguration.

\begin{coro} 
 There is a neat lattice  $\Gamma'\in \mathcal{C}_3$ such that for every $\Gamma''<\Gamma'$, $\bar{X}^{tor}_{\Gamma''}$ has a finite Albanese map. 
 In particular  the universal covering space of $\bar{X}^{tor}_{\Gamma''}$ is Stein.
\end{coro}

\begin{prv}
Corollary \ref{unram} implies that the only curves contracted by the  Albanese map  lie on the boundary. It is thus  enough to replace the $\Gamma'$ in Corollary \ref{unram}  with $\Gamma'\cap \Gamma_{Hirz}$
for which the Albanese map does not contract any boundary curve. 
\end{prv}

\begin{rem} \label{symhirz} \cite [Prop 15.17, P. 155]{DM}
 This configuration of elliptic curves has an order $72$ complex reflection group $H$ of automorphisms which is a $\mu_6$ central extension of $A_4$
 acting as the symmetry group of the confinguration of $4$ points in $\mathbb{P}^1$ given by the points $0, 1, \infty, e^{\frac{2\pi i}{6}}$. 
 $A_4$ acts as an alternating group on the set whose elements are the 4 elliptic curves. 
\end{rem}

\subsubsection{The Orbifold attached to a Normal subgroup of $\Gamma_{Hirz}$} Let us make  Lemma \ref{ramif} more explicit in the present case. 

Let $\Gamma' \triangleleft_{fi} \Gamma_{Hirz}$ be a finite index  normal subgroup. Then $\Gamma'$ is torsion-free and torsion-free at infinity so that the toroidal compactification
$\bar X'\supset   X'=\Gamma' \backslash \chyp $ is a complex projective manifold with an effective $G=\Gamma_{Hirz} \slash \Gamma'$-action which is fixed point free
on $X'$. It is clear that $G\backslash \bar X'= \bar X$ and we denote the corresponding quotient orbifold  by $\mathcal{X}'= [G\backslash \bar X']$
and the corresponding quotient map by $\pi: \bar X' \to \bar X$. 

\begin{lem} \label{lem22}
Let $p'_i \subset E'_i$ be a flag consisting of a point $p_i\in \pi^{-1}(E_i)$ and $E'_i$ the irreducible component of $\pi^{-1}(E_i)=\sum_j F_{ij}$ via $p_i$. 
Consider $$S'_i=Stab_G(p'_i)
< H_i=Stab_G(E_i') < G.$$
Then $S'_i$ is a cyclic central subgroup of $H_i$ of order $d_i=d_i(\Gamma')$, $G_i=H_i\slash S_i'$ acts effectively on $E'_i$ without fixed points, 
$E_i=G_i\backslash E'_i$ and $p^*E_i=\sum_{j\in H_i\backslash G} d_i F_{ij}$. 

Furthermore $\mathcal{X}'$ is equivalent 
to $\mathcal{X}(\bar X, D, d)$ with $d=(d_1, \ldots d_4)$. 
\end{lem}

We shall adopt the notations of Lemma \ref{lem22} in the rest of section \ref{sec2}. 

\begin{lem}
Conversely every proper finite \'etale  mapping $p: \bar Y \to \mathcal{X}(\bar X, D, d)$ comes from a finite index  subgroup $\Gamma'= Im 
(\pi_1(\bar Y \setminus p^{-1} (D)) \to \pi_1(\bar X \setminus D)) $ such that $p^* D= \sum_i d_i Supp(p^{-1} (E_i))$. 
\end{lem}

The support $Supp(D)$  of an effective Cartier divisor $D$ is the sum of its irreducible components with multiplicity one. 

The existence of $\bar Y$ is equivalent to $\mathcal{X}(\bar X, D, d)$ being uniformizable.

\begin{ques}\label{q1}
 For which $d\in \N_{\ge 1}^4$ is $\mathcal{X}(\bar X, D, d)$ uniformizable?  
\end{ques}

\begin{prop}
  $\mathcal{X}(\bar X, D, d)$ is not developpable if $\{d_1, \ldots, d_4\}=\{1,1,1, m \}$ or $\{1,1,m, m' \}$ with $1<m<m'$.
\end{prop}
\begin{prv} 
 In the listed case the exceptional curve in $\bar X$ gives a sub-orbifold equivalent to $\mathbb{P}(m.[0])$ or to $\mathbb{P}(m.[0]+m'[1])$
 which are non developpable, which obstructs the developpability of the ambient orbifold. 
\end{prv}

\begin{lem}\label{q2}
 Assume $d\in \N_{\ge 1}^4$ is such that $\mathcal{X}(\bar X, D, d)$ uniformizable.  Let $Z:=Exc$ be the exceptional curve and $\mathcal{Z}=(\mathbb
 P^1, d_1,d_2,d_3,d_4)$ be 
 the natural suborbifold. $SS(\mathcal{X}(\bar X, D, d)) $ holds iff it holds for $\mathcal{Z}$. 
\end{lem}

\subsubsection{$d_3=d_4=1$}

This corresponds to studying the pair $(\bar X, D'':=D_1+D_2)$. In order to fix notations, we denote by $U''$ the complement of $E_1+E_2$ in $E\times E$ and by 
$X''$ the complement of $D_1+D_2$ in $\bar X$. Plainly $U''=X''\setminus Exc$. 

\begin{lem} The fundamental group of $U''$ is a product of two finite groups on 2 generators
 $\pi_1(U'')=\mathbb F_2(a,b) \times  \mathbb F_2(c,d)$. The fundamental group of $X''$ is the quotient of $\pi_1(U'')$
 by the normal subgroup generated by $[a,b].[c,d]$. 
\end{lem}
\begin{prv}
Elementary calculation. 
\end{prv}

\begin{coro}
 The natural morphism $F_2(a,b) \to \pi_1(X'')$ factors through $H_{\Z}$ the Heisenberg group of unipotent $3\times 3$ upper triangular  matrices with integer coefficients .
\end{coro}
\begin{prv}
 Since $a$ and $b$ commute with $c$ and $d$ the relation $[a,b]=[c,d]^{-1}$ implies that $a$ and $b$ commute with $[a,b]$. But $F_2(a,b)/<<[a,[a,b]], [b,[a,b]]>>\simeq H_{\Z}$.
\end{prv}

We have a non trivial central extension: 
$$
1\to \Z_{[a,b]}=Z(H_{\Z} ) \to H_{\Z}  \to \Z^2_{a,b} \to 1. 
$$

The geometric interpretation is clear. Consider the boundary $B$ of a regular neighborhood of $D_1$ in $\bar X$. 
If the base point is in $B$ $a$ and $b$ can be homotoped in $B\cap U''$. However $B\subset X''$ and $\pi_1(B)=H_{\Z}$ since $D_1^2=-1$. 

\begin{coro}
The fundamental group $\pi_1(X'')$ is the quotient group of $H_{\Z} ^2$ by the diagonal subgroup of its center $\Delta_{\Z}: \Z \to \Z^2$. Thus,  we have  the central extension:
$$ 1 \to \Z \to \pi_1(X'') \to \Z^4=H_1(X'') \to 1, 
$$
$[a,b]$ mapping to a generator of the center and $[c,d]$ to its opposite.
\end{coro}
\begin{prv}
The map $H_{\Z} ^2 \to\pi_1(X'')$ comes from the previous lemma which can  also be applied with $c,d$ and all these groups have naturally isomorphic abelianization. 
\end{prv}

\begin{coro} \label{unif1}
 $\mathcal{X}(\bar X, D, n,n,1,1)$ is uniformizable.
\end{coro}

\begin{prv} The new relations to get $\pi_1(\mathcal{X}(\bar X, D, n,n,1,1))$ are $[a,b]^n=[c,d]^n=1$. 
 It thus suffices to consider the quotient map  $$\pi_1(\mathcal{X}(\bar X, D, n,n,1,1) \to  H^2_{\Z/n\Z}\slash\Delta_{\Z}(\Z/n\Z).$$
\end{prv}

Thanks to lemma \ref{q2} $SSC(\mathcal{X}(\bar X, D, n,n,1,1))$ is trivially true since the fundamental group of $\mathcal Z\simeq \mathbb{P}^1(\sqrt[n]{0+\infty})$ is finite.

\subsubsection{$d_4=1$}

This corresponds to studying the pair $(\bar X, D':=D_1+D_2+D_3)$. In order to fix notations, we denote by $U'$ the complement of $E_1+E_2+E_3$ in $E\times E$ and by 
$X'$ the complement of $D_1+D_2+D_3$ in $\bar X$. Plainly $U'=X'\setminus Exc$. 

\begin{lem}\label{relx'}
 The fundamental group $G'=\pi_1(U')$ has 5 generators $a,b,c,d,e$ and is presented by the relations: 
 $$ a^{-1} ca=c, a^{-1} e a= e, a^{-1} d a = c^{-1} d c e, 
 $$
 $$b^{-1}db=d, b^{-1} e b=e, b^{-1} c b=d^{-1}c d e^{-1}.
 $$
\end{lem}
\begin{prv}
 Omitted. Easy  calculation. 
\end{prv}

It is is a semi-direct product of $F_3(c,d,e)$ by $F_2(a,b)$.

One has to set $\alpha=c^{-1}a$ $\beta=d^{-1}b$ $e=1$ to recover the previous case.

\begin{lem}
 Let $V'\subset U'$ the trace on $U'$ of a regular neighborhood of $Exc$. Then $\pi_1(V')\hookrightarrow \pi_1(U')$ is generated by 
 $a_1,a_2,a_3$ subject to the relations:
 $$ a_3a_2a_1=a_2a_1a_3=a_1a_3a_1$$
 where 
 $$
 a_3=e, \ a_1=c^{-1}d^{-1}cde^{-1}, \ a_2=\alpha^{-1}\beta^{-1}e^{-1}\alpha\beta.
 $$
The center $Z(\pi_1(V'))$ is infinite cyclic generated by the  element $\gamma_Z=a_3a_2a_1=eaba^{-1}b^{-1}.$
\end{lem}

\begin{prv}
 Omitted. In principle the  calculation is easy, but it turned out to be rather messy. 
\end{prv}

\begin{lem}
The fundamental group $G:=\pi_1(X')$ has 5 generators $a,b,c,d,e$ and is presented by the relations \ref{relx'} plus:
$$
eaba^{-1}b^{-1}=1.
$$

The fundamental group $\pi_1(\mathcal{X}(\bar X, D, d_1,d_2,d_3,1))$ has 5 generators $a,b,c,d,e$ and is presented by the relations \ref{relx'} plus:
$$
a_3a_2a_1=a_1^{d_1}= a_2^{d_2}= a_3^{d_3}=1.
$$

The map $\pi_1(\mathcal Z)=F_2(a_1,a_2)/<<a_1^{d_1}, a_2^{d_2},(a_2a_1)^{d_3}>>\to \pi_1(\mathcal{X}(\bar X, D, d_1,d_2,d_3,1))$ maps $a_1,a_2$ to their
above expressions. 
 
\end{lem}

Hence by Lemma \ref{q2} SSC holds true in that  case if and only if we can find a finite index normal subgroup $H$ of $G(d_1,d_2,d_3):=\pi_1(\mathcal{X}(\bar X, D, d_1,d_2,d_3,1))$
and $\eta \in H_1(H, \Q)$ which does not vanish on $<a_1,a_2>\cap H$. The worst possible choice is when  $G(d_1,d_2,d_3)/H$ is abelian. Unfortunately since 
$H_1(G(d_1,d_2,d_3))$ has rank $4$ a lot of the $H$ one gets with
\texttt{LowIndexSubgroups} in GAP or MAGMA have that property.

\begin{lem}\label{unif2}
 $\mathcal{X}(\bar X,D,(n,n,n,1))$ is uniformizable.
\end{lem}

\begin{prv}
 Immediate consequence of Corollary \ref{unif1}. Indeed, thanks to Lemma \ref{morphstack},  we have a map $$\mathcal{X}(\bar X,D,(n,n,n,1)) \to \mathcal{X}(\bar X,D,(n,n,1,1))$$ and a map 
 $\mathcal{X}(\bar X,D,(n,n,n,1)) \to \mathcal{X}(\bar X,D,(1,n,n,1))$ thanks to Remark \ref{symhirz} which gives a group morphism in a finite group
 $$\pi_1(\bar X,D,(n,n,n,1)) 
 \to  (H^2_{\Z/n\Z}\slash\Delta_{\Z}(\Z/n\Z))^2$$
 which is injective on the isotropy groups. 
\end{prv}

Hence $G_n=\pi_1(\mathcal{X}(\bar X, D, n,n,n,1))$ is the fundamental group of a complex projective surface. 
Let us introduce the following quotient of $G_n$:
$$G_n^-:=G_n/<< a^n,b^n, c^n,d^n>>.
$$
where $a,b,c,d$ denote the natural images of the generators of group $G$. 

\begin{lem} The nilpotent group\footnote{We use the following notation for the central series of a group $G$:  
$\gamma_1(G)=G$,  $\gamma_2(G)=[G,G]$, $\gamma_{k+1}(G)=[\gamma_k(G),G]$ } of class $2$ $G_n^-/\gamma_3(G_n^-)$
is isomorphic to $H^2_{\Z / n\Z}$. 
\end{lem}

\begin{prv}
 Using MAGMA \cite{magma}, we get an isomorphism: 
$$
G/\gamma_3(G)\simeq H^2_{\Z}, 
$$
whence the result since we are killing the $n$-th power of the lifted  generators of $H_1(G)=G/\gamma_2(G)$. A file containing the MAGMA code is available on my webpage.
\end{prv}

\begin{prop} $SSC(\mathcal{X}(\bar X,D,(3,3,3,1))$ holds  and the universal covering space is Stein. 
\end{prop}
\begin{prv}

The abelianization $A_3$ of $K_3:=\ker(\phi)$, where we denote by $\phi$ the resulting  epimorphism $\phi : G_3 \to H^2_{\Z / 3\Z}$ is free of rank $10$, thanks to MAGMA. 
MAGMA computes the image of $a_1a_3a_1^{-1}a_3^{-1}$ to be the row vector $(0,-1,-1,1,-1,0,0,0,-1,2)$ in MAGMA's basis of $A_3$.
 
\end{prv}

\begin{coro} For all $k,l, m \in \N^*$,  $SSC(\mathcal{X}(\bar X,D,(3k,3l,3m,1)))$ holds  and the moduli space of the universal covering stack is Stein, hence 
the universal covering space is Stein
provided the orbifold is indeed developable. 
\end{coro}
\begin{prv}
 Use the natural map $\mathcal{X}(\bar X,D,(3k,3l,3m,1))\to \mathcal{X}(\bar X,D,(3,3,3,1))$ given by Lemma \ref{morphstack} to deduce that the Albanese map is virtually finite. 
\end{prv}

\subsubsection{General case} \label{subsec:general}

We have no general results on uniformizability. 
It follows from the above that $\mathcal{X}(\bar X,D,(n,n,n,n))$ is uniformizable for all $n\in \N^*$ , and $SSC(\mathcal{X}(\bar X,D,(3k,3l,3m,3p))$ holds.

\subsection{Some small refined commensurability classes in $\mathcal{C}_3$}

\subsubsection{}
Quite confusingly, the Picard-Eisenstein lattice studied in \cite{Hol2,Hol} is not the same as the one studied in \cite{FP}. Also the relationship with $\Gamma_{Hirz}$  and other 
lattices in $\mathcal{C}_3$ is slightly  involved.  Holzapfel's Picard modular group is (conjugate by a transposition matrix to)
$\Gamma_{H_1, L_{st}}$ for the hermitian form whose matrix is $H_1$  and $L_{st}=\mathcal {O}_3^{\oplus 3}=\Z[\frac{1+i\sqrt{3}}{2}]^{\oplus 3}$ is the standard lattice. 
By definition $H_0=^t{\bar g} H_1 g$ where we have used the notations 
$$ H_1=
\left(
\begin{array}{ccc}
 1&0&0 \\
 0&-1&0 \\
 0&0&1
\end{array}
\right), \ 
g=\left(
\begin{array}{ccc}
 0&1&0 \\
 1&0&1/2 \\
 1&0&-1/2
\end{array}
\right). $$

As $\Q$-algebraic groups 
$PU(H_1)\simeq PU(H_0)$ - more precisely $gPU(H_0)g^{-1}=PU(H_1)$ but the $PU(H_i,\mathcal{O}_3)$  $i=0,1$ are different lattices. 
As we will see, the refined commensurability classes of these two lattices have a very similar behaviour but we did not complete the elementary but lenghty
calculations to check that they are equal. 

\subsubsection{}
The article \cite{Hol2} uses another special elliptic configuration in $E\times E$.  Let $0_E, Q_1, Q_2$ be the fixed points of the automorphism $j: E \to E$
where $\mathcal{O}_3^{\times}$ acts by multiplication on $E=\C\slash \mathcal{O}_3$ and $\mathcal{O}_3=\Z[j]\subset \C$ is the inclusion given by our choice of $i=\sqrt{-1}$
and the formula $j=-\frac{1}{2} + i \frac{\sqrt{3}}{2}$. Then $\{ 0_E, Q_1, Q_2 \}$  defines a group of translations $T:=\mathbb{\Z}\slash 3 \mathbb{Z}$ of $E$. 
Let $T$ acts diagonally on $E\times E$. 
It turns out that there is an isomorphism $T\backslash E\times E \simeq E \times E$ such that the inverse image of the Hirzebruch configuration 
is the union of 6 elliptic curves. This elliptic curve  configuration has  $\{ 0_E\times 0_E ,Q_1 \times Q_1, Q_2 \times Q_2 \}$ as its multiple (quadruple) points. These 6 elliptic
curves $\{ S_k \}_{k=1}^6$ are
the graph of the automorphisms $1,j, j^2$ and the horizontal factors $\{ E\times 0_E, E\times Q_1,  E\times Q_2 \}$. 
Thus there is an index $3$ subgroup $\Gamma_{Holz} \subset \Gamma_{Hirz}$ which can be defined by the Galois correspondance 
$$\Gamma_{Holz}=Im(\pi_1(Bl_{T.(0_E\times 0_E) }(E\times E)\setminus \cup_{k=1}^6 S'_k)\to \pi_1(\bar X\setminus D)=\Gamma_{Hirz})\lhd\Gamma_{Hirz}$$ where 
$S'_k$ is the strict transform of $S_k$ and the notation $(\bar X,D)$ of subsection \ref{sec2} still applies.

One of the facts \cite{Hol2} uses is that we have a composition sequence:

$$ \Gamma_{Holz} \lhd P\tilde\Gamma'_{Holz}:=P(SU(H_1, \mathcal{O}_3)(1-j)) \lhd  PU(H_1, \mathcal{O}_3)$$

with graded quotients $$P\tilde\Gamma'_{Holz}\slash \Gamma_{Holz}=\mu_3 \times \mu_3, \quad PU(H_1, \mathcal{O}_3)\slash P(SU(H_1, \mathcal{O}_3)(1-j))=S_4.$$

Furthermore, thanks to \cite{GerSa},  we can interpret  a crucial ingredient in \cite{Hol} as an equivalence $$ [P\tilde\Gamma'_{Holz}\backslash \chyp] \simeq \mathcal{X}(\mathbb{P}^2\setminus \{ 4 pts\}, {CEVA(2)},(3, 3,3,3,3,3))$$
$$[S_4 \backslash  [P\tilde\Gamma'_{Holz}\backslash \chyp]]\simeq[PU(H_1, \mathcal{O}_3)\backslash\chyp], $$

where $CEVA(2)$ is the complete quadrangle built on the 4 marked points (in general linear position) in $\mathbb{P}^2$ and $S_4<PGL(3,\C)$ permutes these $4$ points. 
\subsubsection{}
Using the equivalence $[\mu_3 \backslash E]\simeq \mathbb{P}^1(\sqrt[3]{0+1+\infty})$ we see easily that the quotient orbifold $[\mu_3\times \mu_3 \backslash Bl_{T.(0_E\times 0_E) }(E\times E)]$ 
is the following orbisurface: the moduli space is $Bl_{(0,0),(1,1), (\infty,\infty)}\mathbb{P}^1 \times \mathbb{P}^1$ and the ramification  has order $3$
on the 9 curves given by the 3 exceptional curves and the strict transforms of the vertical and  horizontal factors through the 3 blown up points. 
There is no ramification over the strict transform of the diagonal. 

In particular $\mathcal{X}_{P\tilde\Gamma'_{Holz}}^{tor}$ is the orbifold whose moduli space is $Bl_{(0,0),(1,1), (\infty,\infty)}\mathbb{P}^1 \times \mathbb{P}^1$
and which ramifies at order 3 on the 6 (-1)-rational curves given by the 3 exceptional curves and the 3 strict transforms of the vertical factors. The $4$ boundary components
are the strict transforms of the diagonal and the 3 strict transforms of the horizontal factors carry multiplicity $1$ since our construction of the 
orbifold toroidal compactifications precisely excludes orbifold behaviour at the general point of the boundary. 

If  $CEVA'(2)$ denotes the strict transform of $CEVA(2)$ in the blow up surface $Bl_{4 pts}(\mathbb{P}^2)$
we see using the familiar contraction of  these 4 disjoint rational curves and  \cite{GerSa} an equivalence: 
$$
\mathcal{X}_{P\tilde\Gamma'_{Holz}}^{tor}\simeq  \mathcal{X}( Bl_{4 pts}(\mathbb{P}^2), {CEVA'(2)}, (3, 3,3,3,3,3))
$$

Using the language of \cite{BHH} we assign  weight 3 to the strict transforms of the lines in $CEVA(2)$ and the weight 1 to the exceptional curves. 
Since the map $\mathcal{X}_{\Gamma_{Holz}}^{tor}\to[\mu_3\times \mu_3 \backslash Bl_{3 pts}(E\times E)]$ is \'etale, 
$\mathcal{X}_{\Gamma_{Holz}} \to \mathcal{X}_{P\tilde\Gamma'_{Holz}}$  is the only non-\'etale map in the orbifold version of  the main diagram in \cite{Hol2}:
$$\xymatrix{
 & &\mathcal{X}_{\Gamma_{Holz}}^{tor}= Bl_{3 \ pts}(E\times E)\ar[dl]_ {ram} \ar[dr]^{et} & \\
 &\mathcal{X}_{P\tilde\Gamma'_{Holz}}^{tor} \ar[dl]_{et}& &\mathcal{X}_{\Gamma_{Hirz}}^{tor}= Bl_{0}(E\times E) \\
\mathcal{X}_{PU(H_1,\mathcal{O}_3)}^{tor}&&&
}$$

In other words, we have in $\mathcal{C}_{3r}$:
$$[PU(H_1,\mathcal{O}_3)]=[P\tilde\Gamma'_{Holz}] \prec [\Gamma_{Holz}]=[\Gamma_{Hirz}]. 
$$
Both classes are small of $\Gamma$-dimensions $1$ and $2$. 

\begin{prop}\label{fib} $(P\tilde\Gamma'_{Holz})^{tor} \cong \pi_1(\mathbb{P}^1(3,3,3))\cong\Z[j]\rtimes \mu_3 $.
$PU(H_1,\mathcal{O}_3)^{tor}$ is virtually abelian of rank 2 sitting in an exact sequence:
$$ \{1\}\to \Z[j]\rtimes \mu_3 \to PU(H_1,\mathcal{O}_3)^{tor} \to S_4 \to \{1\}.
$$
\end{prop}
\begin{prv}
Let us consider the linear system $|I_{4 pts}(2)|$ of conics through the four points. It defines a rational map 
$\phi: \mathbb{P}^2 \dashrightarrow \mathbb{P}^1$ which becomes regular on $Bl_{4pts}(\mathbb{P}^2)$, the exceptional curves are then  isomorphically mapped 
$\phi$ and $CEVA'(2)$ has three connected components which have 2 irreducible  as components and coincide to the three singular fibers of $\phi$. 
The generic fibre is a smooth conic with no deleted points. In other words,  $$\phi: \mathbb{P}^2 \setminus CEVA'(2) \to \mathbb{P}^1 \setminus \{ 3 \ pts \}$$
is a projective smooth conic bundle.  Lemma \ref{morphstack} gives a map:
$$\bar \phi: \mathcal{X}_{P\tilde\Gamma'_{Holz}}^{tor} \to \mathbb{P}^1(3,3,3)
$$
whose general fiber is a smooth rational curve. This gives an isomorphism $$(P\tilde\Gamma'_{Holz})^{tor} \buildrel{\cong }\over\longrightarrow \pi_1(\mathbb{P}^1(3,3,3)).$$ The other statements 
are immediate consequences granted the geometric description above. 
\end{prv}

\begin{coro} When  $\Gamma$ lies in $[PU(H_1,\mathcal{O}_3)]$ or $ [\Gamma_{Hirz}]$, 
$\Gamma^{tor}$  is infinite  virtually abelian and linear.
\end{coro}

\subsubsection{}

In the notations of \cite{Mos},  the group $PU(H_0,\mathcal{O}_3)$ can be described as $\Gamma_{\mu, S_3}$ with $\mu=\frac{1}{6}(2,2,2,1,5)$  \cite{Der}.  Since this ball 
5-uple satisfies INT, the orbifold $[PU(H_0,\mathcal{O}_3)\backslash \chyp]$ can be easily described. 

The moduli space of $[\Gamma_{\mu}\backslash \chyp]$ is $\mathbb{P}^2\setminus \{ P \}$ this point $P$ being
a triple point of $CEVA(2)$. In $[\Gamma_{\mu}\backslash \chyp]$ the  $3$
 lines through $P$ have orbifold weight 3, the three remaining
 lines have orbifold weight 2 and the 3 triple points have a non abelian order $36$ inertia group \cite[pp 392-393]{Ulu}. 
 
Then $\mathcal{X}_{\Gamma_{\mu}}^{tor}$ has $Bl_P(\mathbb{P}^2)$ as its moduli space and the only modification is that we should affect the weight $1$
to the exceptional curve. 

When we mod out be the action of $S_3$ fixing $P$ and permuting the other triple points we observe that there is no ramification
on the exceptional curve. In particular $[\mathcal{X}_{PU(H_0,\mathcal{O}_3)}^{tor}]\simeq S_3 \backslash\mathcal{X}_{\Gamma_{\mu}}^{tor}$. 
In terms of refined coholomogy classes $ [PU(H_0,\mathcal{O}_3)]=[\Gamma_{\mu}]$.

The central projection to $P$ defines a map $\mathcal{X}_{\Gamma_{\mu}}^{tor}\to \mathbb{P}^1(3,3,3)$ whose general fiber $F$ is a $ \mathbb{P}^1(2,2,2)$ an elliptic orbifold
whose fundamental group is the Vierergruppe. This implies that $PU(H_0,\mathcal{O}_3)^{tor}$ is virtually abelian of rank $2$ and that the morphism of Proposition \ref{surj}
has a finite kernel.  Actually,  $\pi_1(F)$  injects thanks to:

\begin{prop}
 There is a (split) exact sequence $$1 \to (P\tilde\Gamma'_{Holz})^{tor} \to \Gamma_{\mu} ^{tor}\to K_4 \to 1$$ and the refined commensurability classes of $PU(H_0,\mathcal{O}_3)$
 and $PU(H_1,\mathcal{O}_3)$ are the same. 
\end{prop}

\begin{prv}
The group $S_4$ acts on the linear system $|I_{4 pts}(2)|$ by the projectivities preserving the 4 points. On the 3 singular members it acts as 
$S_3$ where $S_3=S_4/K_4$ where $K_4$ is the Vierergruppe or Klein group isomorphic to $(\Z/2)^2$. In particular $K_4$ acts as automorphisms of the
map $\phi$ (see also \cite[I.6.2]{Hol}) and its orbifold compactification $\bar \phi$ in Proposition \ref{fib} .  It is then easy to interpret \cite[I.3.6.3]{Hol} 
and see that $[K_4 \backslash \mathcal{X}_{P\tilde\Gamma'_{Holz}}^{tor}]=[\mathcal{X}_{\Gamma_\mu}^{tor}]$. 
\end{prv}

\subsubsection{The universal covering stack attached to $[PU(H_0,\mathcal{O}_3)]$}

\begin{prop} $[\mathcal{X}_{PU(H_0,\mathcal{O}_3)}^{tor}]$ is not developable. 
 \end{prop}
\begin{prv}
 We consider the \'etale map $e:[\mu_3\backslash E]\buildrel{\cong}\over\longrightarrow \mathbb{P}^1(3,3,3)$. Then
 $$ E \times_{e, \bar \phi} \mathcal{X}_{P\tilde\Gamma'_{Holz}}^{tor}\simeq \mathcal{X}_{P\Gamma'_{Holz}}^{tor}$$
 where $P\Gamma'_{Holz}$ is the lattice in $PU(2,1)$ corresponding to the group $\Gamma'$ in the notations of \cite[p. 27]{Hol}. 
 The moduli space of that stack is a surface birational to $E\times\mathbb{P}^1$ with 3 $A_2$ singular points 
 which is actually isomorphic to $\mu_3\times \{1 \} \backslash Bl_{T.0_E\times 0_E} E\times E$. 
 The orbifold structure of $\mathcal{X}_{P\Gamma'_{Holz}}^{tor}$ a $\mu_3$ inertia group at the singular points. 
 So there are orbifold points in the fiber of the map $\mathcal{X}_{P\Gamma'_{Holz}}^{tor} \to E$ which is an isomorphism on $\pi_1$.
 In particular the universal covering stack is equivalent to $\C \times_E \mathcal{X}_{P\Gamma'_{Holz}}^{tor}$ and has infinitely many $\mu_3$ orbifold points. 
\end{prv}

\subsubsection{$SSC(\mathcal{X}(\bar X, D, (n,n,n,n)))$}

The results  in subsection \ref{subsec:general} were slightly unsatisfying but one can settle $SSC$ in the case where the $d_i$ have a common factor:

\begin{prop} If $n\in \N_{\ge 2}$,  $SSC(\mathcal{X}(\bar X, D, (n,n,n,n))$ holds and the universal covering space is a Stein manifold.
 \end{prop}
 
 \begin{prv}
  It is enough to prove this  for the 3-1 \'etale cover $\mathcal{X}'$ given by the blow up at 3 points of $E\times E$ 
  with weights $n$ on
  the strict transform of the Holzapfel configuration of $6$ elliptic curves. Let us take the (\'etale) quotient stack  by the action of $\mu_3\times \mu_3$
we have already encountered. The resulting orbifold $\mathcal{Y}$ can be described in the language of \cite{BHH} as follows:
the moduli space is $\mathbb{P}^2$ blown up at 4 points, the strict transforms of the lines in CEVA(2) have weight 3, 
one exceptional curve has weight $n$,  the 3 other exceptional curves have weight $3n$. 

Let us now consider the linear projection from this blown up $\mathbb{P}^2$ to $\mathbb{P}^1$ with center the point whose exceptional curve has weight $n$. 
It is a regular map which  has 3 special fibers which are isomorphic to a nodal conic, the irreducible constituents carrying weights 3 and $3n$.
Hence there is a map $\mathcal{Y} \to \mathbb{P}^1(3,3,3)$. Composing with the natural \'etale map $\mathcal{X}'\to \mathcal{Y}$, we get a map
$\mathcal{X}'\to \mathbb{P}^1(3,3,3)$ which is not constant on the 3 preimages of $\mathcal{Z}$. Since $\mathbb{P}^1(3,3,3)$ is elliptic and  has virtually abelian 
rank 2 fundamental group
the proposition follows.

 \end{prv}

\section{Concluding remarks}
We conclude by a discussion of some  interesting examples from the litterature. 

In 
\cite{Stov, DS}, bielliptic smooth toroidal compactifications of ball quotients are constructed. They satisfy $L(\_)$ since the fundamental group is virtually abelian hence linear. 
The Shafarevich conjecture is established for surfaces of Kodaira dimension $\le 1$ by \cite{GUR} and their argument gives that the fundamental group is linear. 

More to the point, \cite[Theorem 1.3]{DiCSto} asserts that a smooth toroidal compactification of a ball quotient which is birational to an abelian or a bielliptic surface
is the blow up in finitely many distinct points of the minimal surface. In particular, it has a finite (abelian) cover such that
the  connected components of the Albanese fibres are smooth hence irreducible.  Hence,  one cannot construct a counterexample to $SC(\_)$ by ramifying 
along these connected components as in subsection \ref{sec2}. It is not clear whether $SSC(\_)$ is satisfied.

\noindent {Philippe Eyssidieux}\\
{Universit\'e  de Grenoble-Alpes. Institut Fourier.
100 rue des Maths, BP 74, 38402 Saint Martin d'H\`eres Cedex, France}\\
{philippe.eyssidieux@univ-grenoble-alpes.fr}
{http://www-fourier.ujf-grenoble.fr/$\sim$eyssi/}\\

\end{document}